\documentclass[10pt,reqno]{amsart}

\usepackage{a4wide}
\usepackage{dsfont}                   
\usepackage{amssymb}
\usepackage[all]{xy}
\usepackage{mathrsfs}                 

\newtheorem{proposition}{Proposition}[section]
\newtheorem{theorem}[proposition]{Theorem}

\theoremstyle{remark}
\newtheorem{definition}[proposition]{Definition}
\newtheorem{remark}[proposition]{Remark}
\newtheorem{example}[proposition]{Example}

\newcommand{\cst}{\ifmmode\mathrm{C}^*\else{$\mathrm{C}^*$}\fi}
\newcommand{\st}{\;\vline\;}
\newcommand{\CC}{\mathbb{C}}
\newcommand{\eps}{\varepsilon}
\newcommand{\vt}{\!\vartriangle\!}
\newcommand{\XX}{\mathbb{X}}
\newcommand{\YY}{\mathbb{Y}}
\newcommand{\EE}{\mathbb{E}}
\newcommand{\DD}{\mathbb{D}}
\newcommand{\sS}{\mathbb{S}}
\newcommand{\HH}{\mathbb{H}}
\newcommand{\GG}{\mathbb{G}}
\newcommand{\WW}{\mathbb{W}}
\newcommand{\PP}{\mathbb{P}}
\newcommand{\UU}{\mathbb{U}}
\newcommand{\MM}{\mathbb{M}}
\newcommand{\cH}{\mathscr{H}}
\newcommand{\cF}{\mathcal{F}}
\newcommand{\tens}{\otimes}
\newcommand{\id}{\mathrm{id}}
\newcommand{\comp}{\!\circ\!}
\newcommand{\I}{\mathds{1}}
\newcommand{\Label}[1]{\label{#1}\stepcounter{equation}\tag{\theequation}}
\DeclareMathOperator{\C}{C}
\DeclareMathOperator{\B}{B}
\DeclareMathOperator{\cK}{\mathcal{K}}

\numberwithin{equation}{section}

\begin{document}

\author{Piotr M.~So{\l}tan} \address{Department of Mathematical Methods in Physics, Faculty of Physics, University of Warsaw, Poland}
\email{piotr.soltan@fuw.edu.pl}

\title{Quantum families of maps}

\thanks{Partially supported by National Science Centre (NCN) grant no.~2011/01/B/ST1/05011.}



\begin{abstract}
In this survey article we give basic introduction to the theory of quantum families of maps. We begin with a general look at non-commutative (or ``quantum'') topology. Then we formulate all our results in this language. Existence of quantum families possessing special properties is discussed and we show that these quantum spaces are canonically endowed with quantum semigroup structures. Classical analogy is emphasized at each step and many examples are described.
\end{abstract}


\maketitle

\section{Introduction}

Non-commutative topology is a relatively new branch of mathematics based on principles of noncommutative geometry explained in the seminal book \cite{ncg}. The aim of this paper is to give a survey of the theory of quantum families of maps which constitutes a non-commutative topological attempt to study analogs of spaces of continuous mappings known from general topology. For proofs of the presented results we will refer to the literature.

The paper is organized in the following way: in Section \ref{nctSect} we formulate the Gelfand theorem which lays foundations of non-commutative topology and describe the objects of our investigation, namely \emph{quantum spaces}. Section \ref{qfamSect} deals with the definition of a quantum family of maps and provides a motivation for considering this notion. Furthermore the operation of composition of quantum families of maps is introduced. Then, in Section \ref{allSect} we discuss the case of the non-commutative analog of the space $\C(X,Y)$, where $X$ and $Y$ are topological spaces and the space of continuous mappings $X\to{Y}$ is given the compact-open topology. This object is called the \emph{quantum space of all maps}. We discuss the existence of this quantum space and provide an interesting example. Section \ref{strSect} is devoted to the study of additional structures one obtains from the universal property enjoyed by the quantum space of all maps from a quantum space to itself. The notions of a quantum semigroup and a quantum group are introduced and examples are given. Finally, Section \ref{exmplSect} is devoted to the exposition of various related quantum semigroups including quantum semigroups preserving a state (analog of a probability measure) and quantum commutants. In the last subsection we briefly discuss the quantum semigroup structure on the quantum space of all maps from a finite set to a quantum semigroup (or group).

\section{Non-commutative topology}\label{nctSect}

As mentioned briefly in the introduction, non-commutative topology aims to extend the study of topological spaces to more general objects by viewing the category of spaces as a subcategory of some larger category of objects endowed with extra structure. In our case we will identify the category of compact spaces (and more generally locally compact spaces) with a subcategory of the category of \cst-algebras. The famous theorem of Gelfand discussed in Subsection \ref{imgSub} will provide both the motivation and mathematical background for this point of view. Before this we introduce the notion of a \cst-algebra. In what follows we will always assume that compact and locally compact spaces are Hausdorff.

\subsection{\cst-algebras}

The fundamental notion of non-commutative topology is that of a \cst-algebra.

\begin{definition}
A \cst-algebras is a Banach algebra $A$ over $\CC$ endowed with an anti-linear and anti-multiplicative involutive mapping
\[
A\ni{a}\longmapsto{a^*}\in{A}
\]
such that
\[\Label{cst}
\|a^*a\|=\|a\|^2
\]
for all $a\in{A}$.
\end{definition}

The connection between the norm of a \cst-algebra and the algebraic structure expressed by \eqref{cst} has very far-reaching consequences. The reader is advised to consult classic textbooks on \cst-algebra theory such as e.g.~\cite{dix}.

Examples of \cst-algebras are plentiful:
\begin{itemize}
\item let $\cH$ be a Hilbert space. Then the algebra $\B(\cH)$ of all bounded operators on $\cH$ with the operator norm and hermitian conjugation as involution is a \cst-algebra with unit (the identity operator on $\cH$). This \cst-algebra is non-commutative unless $\cH$ is one dimensional.
\item Let $X$ be a compact space. The space $\C(X)$ of all continuous complex functions on $X$ with pointwise addition and multiplication, the supremum norm and complex conjugation as involution is a commutative \cst-algebra with unit (the constant function $1$).
\end{itemize}

\cst-algebras might not have a unit. For example:
\begin{itemize}
\item let $\cH$ be an infinite-dimensional Hilbert space. Then the subspace $\cK(\cH)$ of $\B(\cH)$ consisting of all \emph{compact} operators $\cH\to\cH$ is a \cst-algebra without unit.
\item Let $X$ be a non-compact locally compact space. Then the space
\[
\C_0(X)=\bigl\{
f\in\C(X)\st\forall\:\eps>0\;\exists\:\text{compact }K\subset{X}\:|f|<\eps\text{ outside }K\bigr\}
\]
called the space of functions \emph{vanishing at $\infty$} is a \cst-algebra with the supremum norm. This \cst-algebra does not have a unit.
\end{itemize}

There are several natural classes of morphisms between \cst-algebras. The most obvious one is that of a all $*$-homomorphisms. More precisely, let $A$ and $B$ be \cst-algebras and let $\Phi\colon{A}\to{B}$. We say that $\Phi$ is a \emph{$*$-homomorphism} if $\Phi$ is linear, multiplicative and $\Phi(a^*)=\Phi(a)^*$ for all $a\in{A}$. If $A$ and $B$ are unital, we often consider only \emph{unital} $*$-homomorphisms, i.e.~those which map the unit of $A$ to the unit of $B$.

\subsection{Gelfand duality}\label{imgSub}

We begin this subsection with a categorical formulation of the famous theorem of Gelfand:

\begin{theorem}[{I.M.~Gelfand \cite{gelfand}}]\label{GT}
The association to each compact space $X$ of the \cst-algebra $\C(X)$ extends to an anti-equivalence of categories of compact spaces with continuous mappings and commutative unital \cst-algebras with unital $*$-homomorphisms.
\end{theorem}

Let us explain more concretely the content of Theorem \ref{GT}. The first thing to note is that the passage from a compact space $X$ to the \cst-algebra $\C(X)$ is a functor, i.e.~whenever we have two compact space $X$ and $Y$ and a continuous map $\phi\colon{X}\to{Y}$ then we have the associated map
\[
\Phi\colon\C(Y)\ni{f}\longmapsto{f\comp\phi}\in\C(X).
\]
This map is easily checked to be a unital $*$-homomorphisms. The second level of Gelfand's theorem is that any unital $*$-homomorphisms from $\C(Y)$ to $C(X)$ is of this form. Moreover any commutative unital \cst-algebra $A$ is isometrically $*$-isomorphic to a \cst-algebra $\C(X)$ for some uniquely determined compact space $X$  (i.e.~there exists a unital $*$-homomorphism $\Phi\colon{A}\to\C(X)$ which is an isometry for the norms on $A$ and $\C(X)$ respectively).

Let us note that mappings on the level of spaces give rise to $*$-homomorphisms of \cst-algebras in the reverse direction. This is the reason why Gelfand's theorem states existence of \emph{anti}-equivalence of categories.

\begin{remark}
Gelfand's theorem can be also formulated for locally compact spaces and not necessarily unital \cst-algebras. However, this version requires a more elaborate notion of a morphism between \cst-algebras.
\end{remark}

Before we proceed let us give an indication of how the passage from a commutative \cst-algebra $A$ to the associated compact space $X$ is achieved.  We consider the set of unital $*$-homomorphisms from $A$ to $\CC$ (the set of complex numbers is canonically a \cst-algebra, in fact $\CC$ is exactly $\C(P)$ where $P$ is a one point space). Such homomorphisms are called \emph{characters} and they form a compact subset $X$ of the dual space $A^*$ of $A$ considered with its weak${}^*$ topology. Then each element $a$ of $A$ is mapped to a continuous function $\hat{a}$ on $X$ defined as
\[
\hat{a}(\lambda)=\lambda(a)
\]
for all $\lambda\in{X}$. The function $\hat{a}$ is called the \emph{Gelfand transform} of $a$ and the mapping $a\mapsto\hat{a}$ from $A$ to $\C(X)$ is also given this name. The proof of Gelfand's theorem consists in part of showing that the Gelfand transform is the desired isometric $*$-isomorphism.

One can consider the set of characters and Gelfand transforms of elements for non-commutative \cst-algebras as well. Of course the set of characters might be empty for some algebras, but in case it is not, we obtain for a \cst-algebras $A$ the universal commutative \cst-algebra $B$ such that any unital $*$-homomorphism from $A$ to a commutative \cst-algebra factorizes through a canonical map from $A$ to $B$. The \cst-algebra $B$ is simply taken to be the image of the Gelfand transform of $A$ sometimes referred to as the \emph{abelianization} of $A$.

Properties of locally compact spaces have sometimes quite obvious reflections on the level of the associated \cst-algebras. We would like to mention two of them:
\begin{itemize}
\item a locally compact space $X$ is compact if and only if the \cst-algebra $\C_0(X)$ is unital,
\item $X$ is finite (has a finite number o points) if and only if $\C(X)$ is finite dimensional (in this case clearly $\C_0(X)=\C(X)$).
\end{itemize}

Lastly let us note that if $X$ and $Y$ are locally compact spaces then the \cst-algebra $\C_0(X\times{Y})$ is canonically isomorphic to the \emph{tensor product} $\C(X)\tens\C(Y)$ of the \cst-algebras $\C(X)$ and $\C(Y)$ (tensor products of \cst-algebras are well described in \cite[Appendix T]{wegge}). In case of commutative \cst-algebras there is only one \cst-algebra completion of the algebraic tensor product, in general we will be using the \emph{minimal} tensor product.

\subsection{Quantum spaces}

Theorem \ref{GT} can be interpreted as saying that the theory of compact spaces is the same as the theory of commutative \cst-algebras. Indeed the two categories are anti-equivalent and any notion pertaining to one class can be expressed using the other. Non-commutative topology is the study of \emph{all} \cst-algebras from the point of view that they are algebras of functions on some spaces. Clearly only commutative ones are \emph{really} algebras of functions, but a lot of interesting properties of commutative algebras have their counterparts in the non-commutative setting and the interpretation that any \cst-algebra is \emph{in some sense} an algebra of functions has far-reaching consequences (the interested reader is referred e.g.~to the paper \cite{ncgnt}).

\begin{definition}
A \emph{quantum space} is an object of the category dual to the category of \cst-algebras.
\end{definition}

We will use a novel notation for quantum spaces. We will use characters like $\XX,\YY,\EE,\DD$ etc. to denote quantum spaces. Each of them corresponds uniquely to a \cst-algebra and the associated (non-commutative) \cst-algebras will be denoted by $\C_0(\XX),\C_0(\YY),\C_0(\EE)$ and $\C_0(\DD)$ respectively. In our language the phrase
\[
\text{``let $\XX$ be a quantum space''}
\]
means exactly the same thing as
\[
\text{``let $\C_0(\XX)$ be a \cst-algebra.''}
\]

In analogy with the last remarks of the previous section we will now make a few definitions:
\begin{itemize}
\item we will say that a quantum space $\XX$ is \emph{compact} if $\C_0(\XX)$ is unital. In this case we will write $\C(\XX)$ for this \cst-algebra.
\item A quantum space $\XX$ will be called \emph{finite} if $\C(\XX)$ is finite-dimensional (in this case $\XX$ is automatically compact.
\end{itemize}

In what follows we will restrict attention solely to compact quantum spaces. In other words all considered \cst-algebras will be unital. Also all $*$-homomorphisms will be unital.

By definition a continuous map from a quantum space $\XX$ to a quantum space $\YY$ is a unital $*$-homomorphism $\Phi\colon\C(\YY)\to\C(\XX)$.

\section{Quantum families of maps}\label{qfamSect}

\subsection{Classical families of maps --- Jackson's theorem}

Consider three sets $A$, $B$ and $C$. It is rather clear that any family of maps $A\to{B}$ indexed by elements of $C$ can be equivalently described as a mapping from $A\times{X}$ to $B$ and vice versa. The next theorem extends this to spaces of continuous mappings of topological spaces.

\begin{theorem}[{J.R.~Jackson \cite{jackson}}]
Let $X$, $Y$ and $E$ be topological spaces such that $X$ is Hausdorff and $E$ is locally compact. For $\psi\in\C(X\times{E},Y)$ define $\sigma(\psi)$ as the mapping from $E$ to $\C(X,Y)$ given by
\[
\bigl[\bigl(\sigma(\psi)\bigr)(e)\bigr](x)=\psi(x,e).
\]
Then $\sigma$ is a homeomorphism of $\C(X\times{E},Y)$ onto $\C\bigl(E,\C(X,Y)\bigr)$ with all spaces of maps topologized by their respective compact-open topologies.
\end{theorem}

In other words a continuous family of maps $X\to{Y}$ parametrized by points of $E$ is encoded in a single map $E\to\C(X,Y)$ and vice versa.

\subsection{Definition of a quantum family of maps}

\begin{definition}\label{qfDef}
Let $\XX$, $\YY$ and $\EE$ be quantum spaces. A continuous \emph{quantum family of maps $\XX\to\YY$} parametrized by $\EE$ is a $*$-homomorphism
\[
\Psi\colon\C(\YY)\longrightarrow\C(\XX)\tens\C(\EE).
\]
\end{definition}

Let us analyze Definition \ref{qfDef} in more detail. The $*$-homomorphism $\Phi$ describes a continuous map from the quantum space corresponding to $\C(\XX)\tens\C(\EE)$ to $\YY$. In view of Jackson's theorem and the fact that tensor product corresponds to Cartesian product of spaces\footnote{This correspondence is not well justified in categorical terms, so we avoid defining the Cartesian product of quantum spaces by taking the tensor product of the respective \cst-algebras.} we can interpret this as a family of maps $\XX\to\YY$ parametrized by the quantum space $\EE$.

Note that the definition of a quantum family of maps is very general. By taking $\C(\EE)$ to be the \cst-algebra $\CC$ we can, for example, treat any $*$-homomorphism as a (trivial) quantum family of maps. In particular it is not difficult to find examples of quantum families of maps.

It is not difficult to prove the following:

\begin{proposition}\label{classqfm}
Let $X,Y$ and $E$ be compact spaces and let $\Phi\colon\C(Y)\to\C(X)\tens\C(E)$ be a $*$-homomorphism. Then there exists a unique family $(\phi_e)_{e\in{E}}$ of maps $X\to{Y}$ parametrized continuously by $E$ such that
\[
\Phi(f)(x,e)=f\bigl(\phi_e(x)\bigr)
\]
for all $x\in{X}$ and $e\in{E}$.
\end{proposition}

Proposition \ref{classqfm} states that when we specify the definition of a quantum family of maps to the ``classical situation'', i.e.~when all considered \cst-algebras are commutative, we obtain the standard notion of a continuous family of maps. More generally, if only $\EE$ is a \emph{classical space} (which means that $\C(\EE)$ is commutative) then there exists a unique family $(\Phi_e)_{e\in{\EE}}$ of $*$-homomorphisms from $\C(\YY)$ to $\C(\XX)$ parametrized continuously (in a natural topology on the set of homomorphisms) by points of $\EE$.

\subsection{Composition of quantum families}

Before we deal with quantum families of maps with additional properties let us define the composition of quantum families of maps. This is a direct generalization of the notion of composition of classical families of maps.

\begin{definition}
Let $\XX_1,\XX_2,\XX_3,\DD_1$ and $\DD_2$ be quantum spaces. Consider families of maps
\[
\begin{split}
\Psi_1\colon\C(\XX_2)&\longrightarrow\C(\XX_1)\tens\C(\DD_1),\\
\Psi_2\colon\C(\XX_3)&\longrightarrow\C(\XX_2)\tens\C(\DD_2).
\end{split}
\]
The \emph{composition} of $\Psi_1$ and $\Psi_2$ is the quantum family of maps
\[
\Psi_1\vt\Psi_2\colon\C(\XX_3)\longrightarrow\C(X_1)\tens\C(\DD_1)\tens\C(\DD_2)
\]
defined by
\[
\Psi_1\vt\Psi_2=(\Psi_1\tens\id)\comp\Psi_2.
\]
\end{definition}

In case the considered quantum families are \emph{classical} (i.e.~the algebras $\C(\DD_1)$ and $\C(\DD_2)$ are commutative) the composition of quantum families corresponds exactly to the family of all compositions of homomorphisms from the families $\Psi_1$ and $\Psi_2$.

Associativity of the operation of composition of quantum families of maps is next to trivial:

\begin{proposition}
Let $\XX_1,\XX_2,\XX_3,\XX_4,\DD_1,\DD_2$ and $\DD_3$ be quantum spaces and let
\[
\begin{split}
\Psi_1\colon\C(\XX_2)&\longrightarrow\C(\XX_1)\tens\C(\DD_1),\\
\Psi_2\colon\C(\XX_3)&\longrightarrow\C(\XX_2)\tens\C(\DD_2),\\
\Psi_3\colon\C(\XX_4)&\longrightarrow\C(\XX_3)\tens\C(\DD_3).
\end{split}
\]
Then $(\Psi_1\vt\Psi_2)\vt\Psi_3=\Psi_1\vt(\Psi_2\vt\Psi_3)$.
\end{proposition}

\section{Quantum families of all maps}\label{allSect}

In this section we consider the non-commutative analog of the space $\C(X,Y)$ with compact-open topology.

\begin{definition}\label{DefAll}
Let $\XX$ and $\YY$, $\EE$ be quantum spaces and let $\Phi\colon\C(\YY)\to\C(\XX)\tens\C(\EE)$ be a quantum family of maps. We say that
\begin{itemize}
\item $\Phi$ is the \emph{quantum family of all maps} from $\XX$ to $\YY$ and
\item $\EE$ is the \emph{quantum space of all maps} from $\XX$ to $\YY$
\end{itemize}
if for any quantum space $\DD$ and any quantum family $\Psi\colon\C(\YY)\to\C(\XX)\tens\C(\DD)$ there exists a unique $\Lambda\colon\C(\EE)\to\C(\DD)$ such that the diagram
\[
\xymatrix{
\C(\YY)\ar[rr]^-{\Phi}\ar@{=}[d]&&\C(\XX)\tens\C(\EE)\ar[d]^{\id\tens\Lambda}\\
\C(\YY)\ar[rr]^-{\Psi}&&\C(\XX)\tens\C(\DD)}
\]
commutes.
\end{definition}

Upon specializing to the classical situation ($\XX,\YY,\EE$ are taken to be classical spaces as well as all possible spaces $\DD$), Definition \ref{DefAll} turns out to be exactly the definition of the space $\C(\XX,\YY)$ with its compact open topology). This result is not entirely trivial.

The natural question which arises in connection with Definition \ref{DefAll} is whether given two quantum spaces $\XX$ and $\YY$ the quantum space of all maps $\XX\to\YY$ exists. The answer is usually negative. The reason for this is that \cst-algebras are only suitable to describe locally compact (quantum or classical) spaces and spaces of continuous mappings with compact open topology rarely are locally compact. Nevertheless there are situations when one can prove existence of the quantum space of all maps. In \cite{pseu} S.L.~Woronowicz stated the following theorem:

\begin{theorem}\label{SLW}
Let $\XX$ and $\YY$ be quantum spaces such that $\C(\XX)$ is finite dimensional and $\C(\YY)$ is finitely generated and unital. Then the quantum space of all maps $\XX\to\YY$ exists. Moreover this quantum space is compact.
\end{theorem}

This means that there always exists the quantum space of all maps from a finite quantum space to one whose associated \cst-algebra is finitely generated and unital. The proof of Theorem \ref{SLW} was given in \cite{qs}.

One way to make sure the hypotheses of Theorem \ref{SLW} are satisfied is to consider a finite quantum space $\MM$ (so $\C(\MM)$ is finite dimensional). Then there always exists the quantum space of all maps $\MM\to\MM$. Contrary to the classical situation, examples of such spaces can be very interesting.

\begin{example}\label{Ex1}
Let $\MM$ be the classical two point space: $\C(\MM)=\CC^2$. Clearly the \emph{classical} space of all maps $\MM\to\MM$ has four points. However the \emph{quantum} space of all maps $\MM\to\MM$ is much more interesting. It turns out to be a quantum space $\EE$ such that
\[
\C(\EE)=\bigl\{f\in\C\bigl([0,1],M_2(\CC)\bigr)\st{f(0),f(1)}\text{ are diagonal}\bigr\}.
\]
The quantum family of all maps $\MM\to\MM$ is $\Phi\colon\CC^2\to\CC^2\tens\C(\EE)$ given by
\[
\Phi\bigl(\left[\begin{smallmatrix}1\\0\end{smallmatrix}\right]\bigr)=
\left[\begin{smallmatrix}1\\0\end{smallmatrix}\right]\tens{P}+
\left[\begin{smallmatrix}0\\1\end{smallmatrix}\right]\tens{Q},
\]
where
\[
P(t)=\begin{bmatrix}0&0\\0&1\end{bmatrix},\quad
Q(t)=\tfrac{1}{2}\begin{bmatrix}
1-\cos{2\pi{t}}&\mathrm{i}\sin{2\pi{t}}\\
-\mathrm{i}\sin{2\pi{t}}&1+\cos{2\pi{t}}
\end{bmatrix}.
\]
The infinite dimensionality of $\C(\EE)$ justifies saying that the quantum space of all maps from a two point space to itself is infinite.

The explanation of these formulas is quite simple: in order to determine a $*$-homomorphism from $\CC^2$ to $\CC^2\tens{A}$, where $A$ is some \cst-algebra, we only need to define its value on $\left[\begin{smallmatrix}1\\0\end{smallmatrix}\right]$. This value can be written as
\[
\left[\begin{smallmatrix}1\\0\end{smallmatrix}\right]\tens{a}+
\left[\begin{smallmatrix}0\\1\end{smallmatrix}\right]\tens{b},
\]
and $a,b\in{A}$ must be projections because $\left[\begin{smallmatrix}1\\0\end{smallmatrix}\right]$ is. The \cst-algebra $\C(\EE)$ described above is simply the universal \cst-algebra generated by two projections without any relations.
\end{example}

\section{Quantum semigroup structure}\label{strSect}

For any set $M$ the set $E$ of all maps from $M$ to $M$ is endowed canonically with the structure of a semigroup. The semigroup multiplication is given by composition of maps. This phenomenon has its non-commutative counterpart.

\begin{theorem}[\cite{qs}]\label{qsemi}
Let $\MM$ be a finite quantum space and let $\EE$ be the quantum space of all maps $\MM\to\MM$. Let
\[
\Phi\colon\C(\MM)\longrightarrow\C(\MM)\tens\C(\EE)
\]
be the quantum family of all maps $\MM\to\MM$. Then there exists a unique
\[
\Delta\colon\C(\EE)\longrightarrow\C(\EE)\tens\C(\EE)
\]
such that $\Phi\vt\Phi=(\id\tens\Delta)\comp\Phi$ i.e.~the diagram
\[
\xymatrix{
\C(\MM)\ar[rr]^-{\Phi}\ar@{=}[d]&&\C(\MM)\tens\C(\EE)\ar[d]^{\id\tens\Delta}\\
\C(\MM)\ar[rr]^-{\Phi\vt\Phi}&&\C(\MM)\tens\C(\EE)\tens\C(\EE)}
\]
is commutative. Moreover the $*$-homomorphism $\Delta$ is coassociative:
\[
(\Delta\tens\id)\comp\Delta=(\id\tens\Delta)\comp\Delta
\]
and there exists a unique character $\epsilon$ of $\C(\EE)$ such that
\[
(\id\tens\epsilon)\comp\Phi=\id
\]
which also satisfies $(\epsilon\tens\id)\comp\Delta=(\id\tens\epsilon)\comp\Delta=\id$.
\end{theorem}

Theorem \ref{qsemi} says that for a finite quantum space $\MM$ the quantum space $\EE$ of all maps $\MM\to\MM$ is canonically endowed with a structure of a \emph{compact quantum semigroup with unit}:

\begin{definition}
Let $\GG$ be a compact quantum space. We say that
\begin{enumerate}
\item $\GG$ is a \emph{compact quantum semigroup} if there exists a coassociative $*$-homomorphism
\[
\Delta\colon\C(\GG)\to\C(\GG)\tens\C(\GG).
\]
We call this morphism the \emph{comultiplication}.
\item $\GG$ is called \emph{compact quantum group} if $\GG$ is a compact quantum semigroup and the sets
\[
\begin{split}
\bigl\{\Delta(a)(\I\tens{b})\st{a,b}\in\C(\GG)\bigr\},\\
\bigl\{(a\tens\I)\Delta(b)\st{a,b}\in\C(\GG)\bigr\}
\end{split}
\]
are linearly dense in $\C(\GG)\tens\C(\GG)$.
\item A compact quantum semigroup $\GG$ has a \emph{unit} if there exists a character $\epsilon\colon\C(\GG)\to\CC$ such that
\[
(\epsilon\tens\id)\comp\Delta=(\id\tens\epsilon)\comp\Delta=\id.
\]
\end{enumerate}
\end{definition}

The definition of a compact quantum group (in the form given above) was introduced by S.L.~Woronowicz in \cite{cqg}. We refer to this paper and to \cite{PodMu} for motivation of this definition and thorough exposition of most fundamental results. Earlier definitions of compact quantum groups were more restrictive (see e.g.~\cite{su2}). These objects are extremely interesting and are currently being studied by many mathematicians. Non-compact quantum groups are fast becoming the avantgarde of current research in harmonic analysis (\cite{timm,mmu}).

It is important to note that unless $\MM$ is a (classical) one point space, the quantum semigroup of all maps $\MM\to\MM$ is never a compact quantum group with its canonical comultiplication (\cite{kom}). In the situation of Example \ref{Ex1} there exists a different comultiplication $\Delta'$ on $\C(\EE)$ with which $\EE$ is a compact quantum group:
\[
\begin{split}
\Delta'(P)&=(P-\I)\tens{P}+\I\tens\I+P\tens(P-\I),\\
\Delta'(Q)&=(Q-\I)\tens{Q}+\I\tens\I+Q\tens(Q-\I).
\end{split}
\]
Existence of a quantum group structure on various quantum spaces is discussed at length in \cite{non,nonBC}.

We have so far established that for a finite quantum space $\MM$ there exists a canonical quantum semigroup structure on the quantum space of all maps $\MM\to\MM$. This structure can be further studied and many of its aspects are explained in the paper \cite{qs}. Let us note the following:

\begin{proposition}[\cite{qs}]
Let $\MM$ be a finite quantum space and let $\EE$ be the quantum space of all maps $\MM\to\MM$. Then the Gelfand transform of $\C(\EE)$ maps $\CC(\EE)$ onto the \cst-algebra of all continuous functions on the compact space of all unital $*$-homomorphisms $\C(\MM)\to\C(\MM)$.
\end{proposition}

\section{Further constructions and examples}\label{exmplSect}

Having established that the quantum space of all maps from a finite quantum space to itself is a quantum semigroup we now proceed to study its subsemigroups.

\subsection{Quantum families preserving a state}

If $M$ is a finite set and $\mu$ is a measure on $M$ then the set of all maps $M\to{M}$ preserving the measure $\mu$ is a subsemigroup of the semigroup of all maps $M\to{M}$. An analogous result is true for finite quantum spaces.

In order to see this we need to introduce a non-commutative analog of a measure on a finite space. For simplicity we will restrict attention to probability measures.

\begin{definition}
Let $A$ be a \cst-algebra and let $\omega$ be a continuous linear functional on $A$. We say that $\omega$ is a \emph{state} if $\|\omega\|=1$ and for any $a\in{A}$ we have $\omega(a^*a)\geq{0}$.
\end{definition}

States on commutative \cst-algebras (i.e.~algebras of functions on compact spaces) correspond to integration with respect to probability measures. We now introduce the notion of a quantum family preserving a given state.

\begin{definition}
Let $\MM$ be a finite quantum space and let $\omega$ be a state on $\C(\MM)$. Let $\DD$ be another quantum space and let $\Psi\colon\C(\MM)\to\C(\MM)\tens\C(\DD)$ be a quantum family of maps $\MM\to\MM$. We say that $\Psi$ \emph{preserves} $\omega$ if
\[
(\omega\tens\id)\bigl(\Psi(x)\bigr)=\omega(x)\I,
\]
for all $x\in\C(\MM)$.
\end{definition}

\begin{theorem}[\cite{qs}]
Let $\MM$ be a finite quantum space and $\omega$ a state on $\C(\MM)$. Then
\begin{enumerate}
\item
there exists a unique quantum family
\[
\Phi_\omega\colon\C(\MM)\longrightarrow\C(\MM)\tens\C(\WW)
\]
such that for any quantum family
\[
\Psi\colon\C(\MM)\longrightarrow\C(\MM)\tens\C(\DD)
\]
preserving $\omega$ there exists a unique $\Lambda\colon\C(\WW)\to\C(\DD)$ such that
\[
\xymatrix{
\C(\MM)\ar[rr]^-{\Phi_\omega}\ar@{=}[d]&&\C(\MM)\tens\C(\WW)\ar[d]^{\id\tens\Lambda}\\
\C(\MM)\ar[rr]^-{\Psi}&&\C(\MM)\tens\C(\DD)}
\]
\item $\Phi_\omega$ preserves $\omega$,
\item $\WW$ is a compact quantum semigroup with unit and if $\EE$ is the quantum semigroup of all maps $\MM\to\MM$ then the canonical map $\Lambda\colon\C(\EE)\to\C(\WW)$ intertwines the comultiplications of $\C(\EE)$ and $\C(\WW)$.
\end{enumerate}
\end{theorem}

Quantum semigroups preserving a state are also very interesting. Unlike Example \ref{Ex1} the following one is purely non-commutative.

\begin{example}\label{Ex2}
Let $\C(\MM)=M_2(\CC)$ and choose a parameter $q\in]0,1]$. Let $\omega_q$ be the state on $\C(\MM)$ defined by
\[
\omega_q\left(\left[\begin{smallmatrix}a&b\\c&d\end{smallmatrix}\right]\right)=\tfrac{a+q^2d}{1+q^2}
\]
The quantum semigroup $\WW$ of all maps $\MM\to\MM$ preserving $\omega_q$ can be described as follows: $\C(\WW)$ is generated by three elements $\beta$, $\gamma$ and $\delta$ such that
\[
\begin{aligned}
q^4\delta^*\delta+\gamma^*\gamma+q^4\delta\delta^*+\beta\beta^*&=\I,
&\beta\gamma&=-q^4\delta^2,\\
\beta^*\beta+\delta^*\delta+\gamma\gamma^*+\delta\delta^*&=\I,
&\gamma\beta&=-\delta^2,\\
\gamma^*\delta-q^2\delta^*\beta+\beta\delta^*-q^2\delta\gamma^*&=0,
&\beta\delta&=q^2\delta\beta,\\
q^4\delta\delta^*+\beta\beta^*+q^2\gamma\gamma^*+q^2\delta\delta^*&=\I,&\delta\gamma&=q^2\gamma\delta\\
q^4\delta^*\delta+\gamma^*\gamma+q^2\beta^*\beta+q^2\delta^*\delta&=q^2\I.
\end{aligned}
\]
The comultiplication $\Delta\colon\C(\WW)\to\C(\WW)\tens\C(\WW)$ is given by
\[
\begin{split}
\Delta(\beta)&=q^4\delta\gamma^*\tens\delta-q^2\beta\delta^*\tens\delta+\beta\tens\beta
+\gamma^*\tens\gamma-q^2\delta^*\beta\tens\delta+\gamma^*\delta\tens\delta,\\
\Delta(\gamma)&=q^4\gamma\delta^*\tens\delta-q^2\delta\beta^*\tens\delta+\gamma\tens\beta
+\beta^*\tens\gamma-q^2\beta^*\delta\tens\delta+\delta^*\gamma\tens\delta,\\
\Delta(\delta)&=-q^2\gamma^*\gamma\tens\delta-q^2\delta\delta^*\tens\delta+\delta\tens\beta
+\delta^*\tens\gamma+\beta^*\beta\tens\delta+\delta^*\delta\tens\delta.
\end{split}
\]
The counit $\epsilon$ maps $\gamma$ and $\delta$ to $0$ and $\beta$ to $1$.
\end{example}

The quantum semigroup preserving the state $\omega_q$ described in Example \ref{Ex2} contains the largest (i.e.~the relevant \cst-algebra is the quotient of $\C(\EE)$ by the smallest ideal) quantum \emph{group} which preserves $\omega_q$. This quantum group has been identified in \cite{so3} as the \emph{quantum $\mathrm{SO}(3)$ group} defined by P.~Podle\'s in \cite{podPHD} (see also \cite{podles}). The original definition of this quantum group used representation theory of the quantum $\mathrm{SU}(2)$ group (\cite{su2}). This was not a very satisfactory definition. Another attempt at understanding this quantum group was made in \cite{podles} where it was fount that the \cst-algebra of functions on the quantum $\mathrm{SO}(3)$ group is generated by five elements $A,C,G,K,L$ satisfying
\[
\begin{aligned}
L^*L&=(\I-K)(\I-q^{-2}K),&{\qquad}AK&=q^2KA,\\
LL^*&=(\I-q^2K)(\I-q^4K),&{\qquad}CK&=q^2KC,\\
G^*G&=GG^*,&{\qquad}LG&=q^4GL,\\
K^2&=G^*G,&{\qquad}LA&=q^2AL,\\
A^*A&=K-K^2,&{\qquad}AG&=q^2GA,\\
AA^*&=q^2K-q^4K^2,&{\qquad}AC&=CA,\\
C^*C&=K-K^2,&{\qquad}LG^*&=q^4G^*L,\\
CC^*&=q^2K-q^4K^2,&{\qquad}A^2&=q^{-1}LG,\\
LK&=q^4KL,&{\qquad}A^*L&=q^{-1}(\I-K)C,\\
GK&=KG,&{\qquad}K^*&=K.
\end{aligned}
\]
The map from $\C(\EE)$ (of Example \ref{Ex2}) onto the \cst-algebra generated by $A,C,G,K$ and $L$ is given by
\[
\beta\longmapsto{L},\quad\gamma\longmapsto{-qG},\quad\delta\longmapsto{q^{-1}A}.
\]

\subsection{Quantum commutants}

Consider the following situation: we are given a finite set $M$ and a family $\cF$ of maps $M\to{M}$. Then the set of all maps $M\to{M}$ commuting with elements of $\cF$ is a semigroup under composition of maps. This phenomenon also has its non-commutative analog. In order to describe it we first introduce the notion of commutation of quantum families of maps.

\begin{definition}
Let $\MM$ be a finite quantum space and let
\[
\Psi_1\colon\C(\MM)\to\C(\MM)\tens\C(\DD_1),\qquad
\Psi_2\colon\C(\MM)\to\C(\MM)\tens\C(\DD_2)
\]
be two quantum families of maps. We say that $\Psi_1$ and $\Psi_2$ \emph{commute} if
\[
(\id\tens\sigma)\comp(\Psi_1\vt\Psi_2)=\Psi_2\vt\Psi_1,
\]
where $\sigma$ is the flip
\[
\C(\DD_1)\tens\C(\DD_2)\ni{x\tens{y}}\longmapsto{y\tens{x}}\in\C(\DD_2)\tens\C(\DD_1).
\]
\end{definition}

Given a quantum family of maps of a finite quantum space to itself there always exists the quantum space of all maps commuting with this family. It is canonically a compact quantum semigroup:

\begin{theorem}[\cite{qs}]
Let $\MM$ be a finite quantum space and $\Psi\colon\C(\MM)\to\C(\MM)\tens\C(\DD)$ a quantum family of maps $\MM\to\MM$. Then
\begin{enumerate}
\item there exists a unique quantum family
\[
\Phi_\Psi\colon\C(\MM)\longrightarrow\C(\MM)\tens\C(\UU)
\]
such that for any quantum family $\Theta\colon\C(\MM)\to\C(\MM)\tens\C(\PP)$
commuting with $\Psi$ there exists a unique $\Lambda\colon\C(\UU)\to\C(\PP)$ such that
\[
\xymatrix{
\C(\MM)\ar[rr]^-{\Phi_\Psi}\ar@{=}[d]&&\C(\MM)\tens\C(\UU)\ar[d]^{\id\tens\Lambda}\\
\C(\MM)\ar[rr]^-{\Theta}&&\C(\MM)\tens\C(\PP)}
\]
\item $\Phi_\Psi$ commutes with $\Psi$,
\item $\UU$ is a compact quantum semigroup with unit and if $\EE$ is the quantum semigroup of all maps $\MM\to\MM$ then the canonical map $\Lambda\colon\C(\EE)\to\C(\UU)$ intertwines the comultiplications of $\C(\EE)$ and $\C(\UU)$.
\end{enumerate}
\end{theorem}

Of course any classical family of morphisms $\C(\MM)\to\C(\MM)$ can be interpreted as a quantum family (parametrized by a classical space) we can consider quantum commutants also of classical families --- even consisting of a single morphisms.

\begin{example}
Let, as in Example \ref{Ex2}, $\MM$ be a quantum space such that $\C(\MM)=M_2(\CC)$. Let $\UU$ be the commutant of the (classical) family of maps $\MM\to\MM$ consisting of the single automorphism of $\C(\MM)$ defined by
\[
\psi\colon\left[\begin{smallmatrix}a&b\\c&d\end{smallmatrix}\right]\longmapsto
\left[\begin{smallmatrix}d&c\\b&a\end{smallmatrix}\right].
\]
This family is described in our language by the $*$-homomorphism
\[
\Psi\colon\C(\MM)\longrightarrow\C(\MM)\tens\CC
\]
given by $\Psi(m)=\psi(m)\tens{1}$. Let $\UU$ be the quantum commutant of $\Psi$. The \cst-algebra $\C(\UU)$ is generated by $\alpha,\beta$ and $\gamma$ with
\[
\beta=\beta^*,\quad\gamma=\gamma^*
\]
and
\[
\begin{aligned}
\alpha^*\alpha+\gamma^2+\alpha\alpha^*+\beta^2&=\I,
&\alpha^2+\beta\gamma&=0,\\
\alpha^*\beta+\gamma\alpha^*+\alpha\gamma+\beta\alpha&=0,
&\alpha\beta+\beta\alpha^*&=0,\\
\gamma\alpha+\alpha^*\gamma&=0.
\end{aligned}
\]
The comultiplication acts on generators as
\[
\begin{split}
\Delta(\alpha)&=\I\tens\alpha+(\alpha^*\alpha+\gamma^2)\tens(\alpha^*-\alpha)+
\alpha\tens\beta+\alpha^*\tens\gamma,\\
\Delta(\beta)&=(\alpha\gamma+\beta\alpha)\tens(\alpha-\alpha^*)+\beta\tens\beta+\gamma\tens\gamma,\\
\Delta(\gamma)&=(\beta\alpha+\alpha\gamma)\tens(\alpha^*-\alpha)+\gamma\tens\beta+\beta\tens\gamma,
\end{split}
\]
It can be shown (\cite{kom}) that $\UU$ is not a compact quantum group (with this $\Delta$).
\end{example}

\subsection{Quantum maps into a quantum semigroup}

In a preprint \cite{mms} yet another construction of quantum semigroup structures on quantum families of maps was proposed. It is based on the simple fact that if $K$ is a set  and $S$ is a semigroup then the set of all maps $K\to{S}$ is a semigroup under pointwise multiplication. This simple construction was extended to the quantum family of all maps from a \emph{classical} finite space into a quantum semigroup.

\begin{theorem}\label{qms}
Let $K$ be a finite set and let $\sS$ be a quantum semigroup with comultiplication $\Delta_\sS\colon\C(\sS)\to\C(\sS)\tens\C(\sS)$. Let $\HH$ be the quantum space of all maps $K\to\sS$ and let
\[
\Phi\colon\C(\sS)\to\C(K)\tens\C(\HH)
\]
be the quantum family of all these maps. Then $\HH$ admits a structure of a compact quantum semigroup described by $\Delta_\HH\colon\C(\HH)\to\C(\HH)\tens\C(\HH)$ defined by the diagram
\[
\xymatrix
{
\C(\sS)\ar[rr]^-\Phi\ar[d]_{\Delta_\sS}&&\C(K)\tens\C(\HH)\ar[d]^{\id\tens\Delta_\HH}\\
\C(\sS)\tens{\C(\sS)}\ar[d]_{\Phi\tens\Phi}&&\C(K)\tens\C(\HH)\tens\C(\HH)\\
\C(K)\tens\C(\HH)\tens\C(K)\tens\C(\HH)\ar[rr]_{\id\tens\sigma\tens\id}
&&\C(K)\tens\C(K)\tens\C(\HH)\tens\C(K)\ar[u]_{\mu\tens\id\tens\id}
}
\]
where $\sigma$ is the flip $\C(\HH)\tens\C(K)\to\C(K)\tens\C(\HH)$ and $\mu\colon\C(K)\tens\C(K)\to\C(K)$ is the multiplication map.
\end{theorem}

In \cite{qmqs} a more detailed analysis of the quantum semigroup structure described in Theorem \ref{qms} was carried out. In particular it was shown that $\HH$ is a quantum group if and only if $\sS$ is. Also an example was presented where $K$ is replaced by a non-classical space. It turns out that in this case $\HH$ no longer needs to be a quantum group even when $\sS$ is a classical finite group.

\end{document}